\documentclass[11pt,a4paper]{article}
\usepackage[utf8]{inputenc}

\usepackage{inputenc, amssymb, amsmath}
\usepackage{amsthm}

\usepackage{graphics}

\newtheorem{thm}{Theorem}
\newtheorem*{thrm}{Theorem}

\newtheorem{rmk}{Note}
\newtheorem{ex}{Example}

\DeclareMathOperator{\Cent}{Cent}

\renewcommand{\phi}{\varphi}

\newcommand{\al}{\alpha}

\newcommand{\si}{\sigma}

\newcommand{\N}{\mathbb{N}}
\newcommand{\Z}{\mathbb{Z}}

\renewcommand{\leq}{\leqslant}
\renewcommand{\geq}{\geqslant}

\begin{document}
 
\begin{center} \large Kanunnikov A. L.\footnote{The research was supported by Russian Science Foundation,
Grant 16-11-10013.} \\ \normalsize Lomonosov Moscow State University, the Faculty of Mechanics and Mathematics 
\end{center}

\begin{center}
 \bf \large Criteria on group $G$ for Goldie theorems to be true for all $G$-graded rings
\end{center}

\vskip 5 pt

{\small {\bf Abstract.} We present two criteria for a group $G$ to satisfy the following statements: any $G$-graded gr-prime (gr-semiprime) right 
gr-Goldie ring admits a~gr-semisimple graded right classical quotient ring. The criterion for gr-semiprime rings is that the group $G$ is periodic. 
Actually, the sufficiency of  periodicity was proved by the author in 2011 and the necessity of it follows from 
the well-known counterexample (1979). The main result of the paper concerns the gr-prime case. In this case, Goodearl and Stafford proved the graded 
version of Goldie's Theorem for rings graded by an Abelian group (2000). Developing their idea we extend the class of groups to the following: 
$(\forall g,h\in G)(\exists n\in \N)(gh^n=h^ng)$. Moreover, for any group $G$ outside this 
class, we construct a~counterexample, precisely, a~$G$-graded ring $R$ such that $Q^{gr}_{cl}(R)=R$ is gr-Noetherian but not gr-Artinian.}

\vskip 5 pt


\vskip 5 pt

{\bf Notation.} 
\begin{itemize}
\item $G$ is a group with neutral element $e$.
\item $O(g)$ is the order of element $g\in G$.
\item $R=\bigoplus_{g\in G} R_g$ is an associative ring graded by $G$; we suppose that $1\in R$.
\item $h(R)=\bigcup_{g\in G} R_g$ is the set of homogeneous elements of $R$; $\deg r= g \Leftrightarrow r\in R_g \setminus 0$.
\item gr-P, where P is some ring property, is a standard graded analogue of P. For example, a graded ring $R$ is called gr-semiprime if it does not 
contain nonzero graded nilpotent ideals; $R$ is called a right gr-Goldie ring if it does not contain infinite direct sums of graded right ideals 
and also satisfies the maximal condition for graded right annihilators.
\item $Q^{gr}_{cl}(R)=RS^{-1}$ is the graded classical right quotient ring of $R$, $S$ is a set of homogeneous regular elements of $R$. 
As well known, $Q^{gr}_{cl}(R)$ exists iff Ore's conditions hold for $h(R)$.
\item $Q^{gr}(R)=\varinjlim\limits_{I\in \mathcal{I}} \mathrm{HOM}_R(I,R)$ is the maximal graded right quotient  ring of $R$, $\mathcal{I}$ is a set 
(graded topology) of all gr-dense right ideals of $R$.
\item $M_n(R)(g_1,\ldots, g_n)$, where $n\in \N$, $g_1,\ldots, g_n\in G$, is a matrix ring with the following $G$-grading: 
$M_n(R)(g_1,\ldots, g_n)_h=(R_{g_i^{-1}hg_j})_{ij}$.
\end{itemize}

\vskip 5 pt

\section{Introduction}

The paper is devoted to the graded version of the following well-known theorem.

\begin{thrm}[Goldie, 1960, \cite{G}] The following conditions on ring $R$ are equivalent:

(1) $R$ is a semiprime right Goldie ring;

(2) a right ideal of $R$ is essential iff it contains a regular element;

(3) the right classical quotient ring $Q_{cl}(R)$ of $R$ exists and is semisimple.

Under these conditions, $Q_{cl}(R)$ is simple iff $R$ is prime. % and, besides, the maximal quotient ring $Q(R)$ of $R$ coincides with $Q_{cl}(R)$.
\end{thrm}

Condition (2) in the proof of Goldie's Theorem plays a key role. Johnson \cite{J} established
that this condition together with right finite dimensionality of $R$ gave a sufficient condition
for the coincidence of the classical $Q_{cl}(R)$ and the maximal $Q(R)$ right quotient rings of $R$.

The analogs of this theorem for rings graded by group have been studied by many authors since the late 1970's. 
Review the most significant results.

\begin{itemize}
 \item \cite{N2} There exists a $\Z$-graded gr-semiprime gr-Goldie ring $R$ such that $Q^{gr}_{cl}(R)=R$ is not gr-Artinian:
$$R=k[x,y]/(xy) \;(k \text{ is a field}); \quad R_n=\begin{cases} kx^n & \text{ if } n>0, \\ k & \text{ if } n=0, \\ ky^n & \text{ if } n<0.
                                                     \end{cases}$$
Thus, the graded analogue of implication (1) $\Rightarrow$ (3) does not hold without extra conditions.

\item \cite{Kan2} The graded analogs of $(1)\Leftarrow (2)\Leftrightarrow (3)$ hold similar to the ungraded case. 
The important difference is that the quotient rings $Q^{gr}_{cl}(R)$ and $Q^{gr}(R)$ of gr-semiprime right gr-Goldie ring 
$R$ can be not equal, but $Q^{gr}(R)$ is gr-semisimple anyway. In fact, $Q^{gr}(R)=Q^{gr}_{cl}(R)$ iff $Q^{gr}_{cl}(R)$ exists and is gr-semisimple. 
For ring $R$ from the example considered above, the $Q^{gr}(R)$ is a direct sum of two graded fields, precisely,
$$Q^{gr}(R)=k[x,x^{-1}]\oplus k[y,y^{-1}], \quad Q^{gr}(R)_n=\begin{cases} kx^n\oplus ky^{-n} & \text{ if } n>0, \\ k\oplus k & \text{ if } n=0, \\ 
kx^{-n}\oplus ky^n & \text{ if } n<0. \end{cases}$$
Note that $R=\{(f(x),g(y))\in k[x]\oplus k[y]\subset Q^{gr}(R)\mid f(0)=g(0)\}$ and $x^{-1}$ is an equivalence class of the following homomorphism of $R$-modules: 
$(x,y)\to R$, $x\mapsto 1,y\mapsto 0$.

\item \cite{GS} Let $R$ be a gr-prime gr-Goldie ring graded by an Abelian group $G$. Then there exists the gr-simple gr-Artinian ring $Q^{gr}_{cl}(R)$.
 See also \cite{GS} and the bibliography stated there for applications to quantized coordinate rings.

\item \cite{Kan} Let $R$ be a right graded gr-Goldie ring satisfying the following condition: 
$$R \text{ is right } e\text{-faithful, i.~e.,} \,(\forall g\in G)(\forall r\in R_g\setminus 0)(\exists r'\in R_{g^{-1}})(rr'\ne 0), \text{ and } R_e 
\text{ is semiprime.} \eqno(*)$$
Then the following is true: 
$$\begin{matrix}
   \text{there exists the gr-semisimple ring $Q^{gr}_{cl}(R)=RS_e^{-1}$, where } \\
\text{$S_e=S\cap R_e$ equals the set of regular elements of $R_e$,} \\ \text{$R_e$ is a semiprime right Goldie ring and 
$Q_{cl}(R_e)=Q^{gr}_{cl}(R)_e=R_eS_e^{-1}$.}
  \end{matrix}
\eqno (**)$$

Note that in \cite{NN} there was proved a similar theorem under the following assumptions: $R$ is right $e$-faithful and $R_e$ is a semiprime right 
Goldie ring.

\item \cite{Kan} Condition $(*)$ holds for gr-semiprime $G$-graded right gr-Goldie ring $R$ in both cases: either the support 
$\mathrm{Supp}(R)=\{g\in G\mid R_g\ne 0\}$ of $R$ is finite, or the group $G$ is periodic.
\end{itemize}

In this paper we prove the following two criteria on group $G$ for Goldie theorem to be true for all 
$G$-graded gr-semiprime and, respectively, gr-prime right gr-Goldie rings.

\begin{thm} \label{T1} The following conditions on group $G$ are equivalent:

(1) for any $G$-graded gr-semiprime right gr-Goldie ring $R$, the ring $Q^{gr}_{cl}(R)$ exists and is gr-semisimple;

(2) $G$ is periodic.

Moreover, any gr-semiprime right gr-Goldie ring graded by a periodic group satisfies $(**)$.
\end{thm}

\begin{thm} \label{T2} The following conditions on group $G$ are equivalent:

(1) for any $G$-graded gr-prime right gr-Goldie ring $R$, the ring $Q^{gr}_{cl}(R)$ exists and is gr-simple;

(2) $(\forall g,h\in G)(\exists n\in \N)\;(gh^n=h^ng)$;

(2)$'$ $(\forall g,h\in G)(\exists m,n\in \N)\;(gh^m=h^ng)$.
\end{thm}

Let us discuss about the class of groups from condition (2) and its significant subclasses. 

\begin{rmk}
Before the author proved Theorem \ref{T2}, his student Dmitry Bazhenov had proved \cite{DB} that condition~(1) holds for any group with finite commutator subgroup 
$G'$. As A. U. Olshanskii noted, such groups lie in the variety  satisfying the identity $gh^n=h^ng$ for some $n$. Indeed, if $|G'|=k\in \N$ then 
$$|G:\Cent G'|=\left| G:\bigcap_{g\in G'}\Cent(g)\right| \leqslant \prod_{g\in G'} |G:\Cent(g)|= \prod_{e\ne g\in G'} |\{xgx^{-1}\mid x\in G\}|\leqslant k^{k-1}.$$
Now, let $g,h\in G$. Since $h^{k^{k-1}}\in \Cent{G'}$ and $x^k=e$ for all $x\in G'$, we have $[g,h^{k^k}]=[g,h^{k^{k-1}}]^k=e$
because $[x,y^m]=[x,y]^m$ whenever $y$ commutes with $[x,y]$. So $gh^n=h^ng$ for $n=k^k$.

As corollary, if $|G'|<\infty$ then $G/Z(G)$ is periodic.

The considered identities have a significant place in the theory of group varieties, see \cite[\S 31]{Ol}. 
We also note that free groups of this variety for sufficiently large odd $n$ are torsion-free \cite[Theorem 31.2]{Ol}. 

\end{rmk}

\begin{rmk}
The condition ``$G/Z(G)$ is periodic'' is stronger than condition (2) and is obtained from it by permutation of two quantifiers. 
In fact, condition (2) holds iff any finitely generated subgroup of $G$ admits a~periodic quotient group modulo center. In particular, 
(2) $\Leftrightarrow G/Z(G)$ is periodic whenever $G$ is finitely generated. For arbitrary group, this is not true as the following example shows.
\end{rmk}

\begin{ex}[A. A. Klyachko] Consider the following subgroup of the direct product $\prod_{n\in \N} D_{2n+1}$ of dihedral groups:
$$\mathcal{D}=\{(g_1,g_2,\ldots)\mid \text{almost all } g_i\text{'s are rotations\,}\}.$$ 
It is clear that $Z(\mathcal{D})=\{e\}$ and that the order of the element $\left(R_{\frac{2\pi}{2n+1}}\right)_{n\in \N}$ is infinite. 
($R_\al$ is a rotation through angle $\al$.) 
At the same time, the group $\mathcal{D}$ satisfies condition (2): for given $g,h\in \mathcal{D}$, one can take $n=2(2j_1+1)\ldots (2j_k+1)$, 
where $j_1,\ldots, j_k$ are numbers of all reflections of the sequence $g$.
\end{ex}

See all classes of groups considered above on the following diagram.

\begin{center}
 \includegraphics{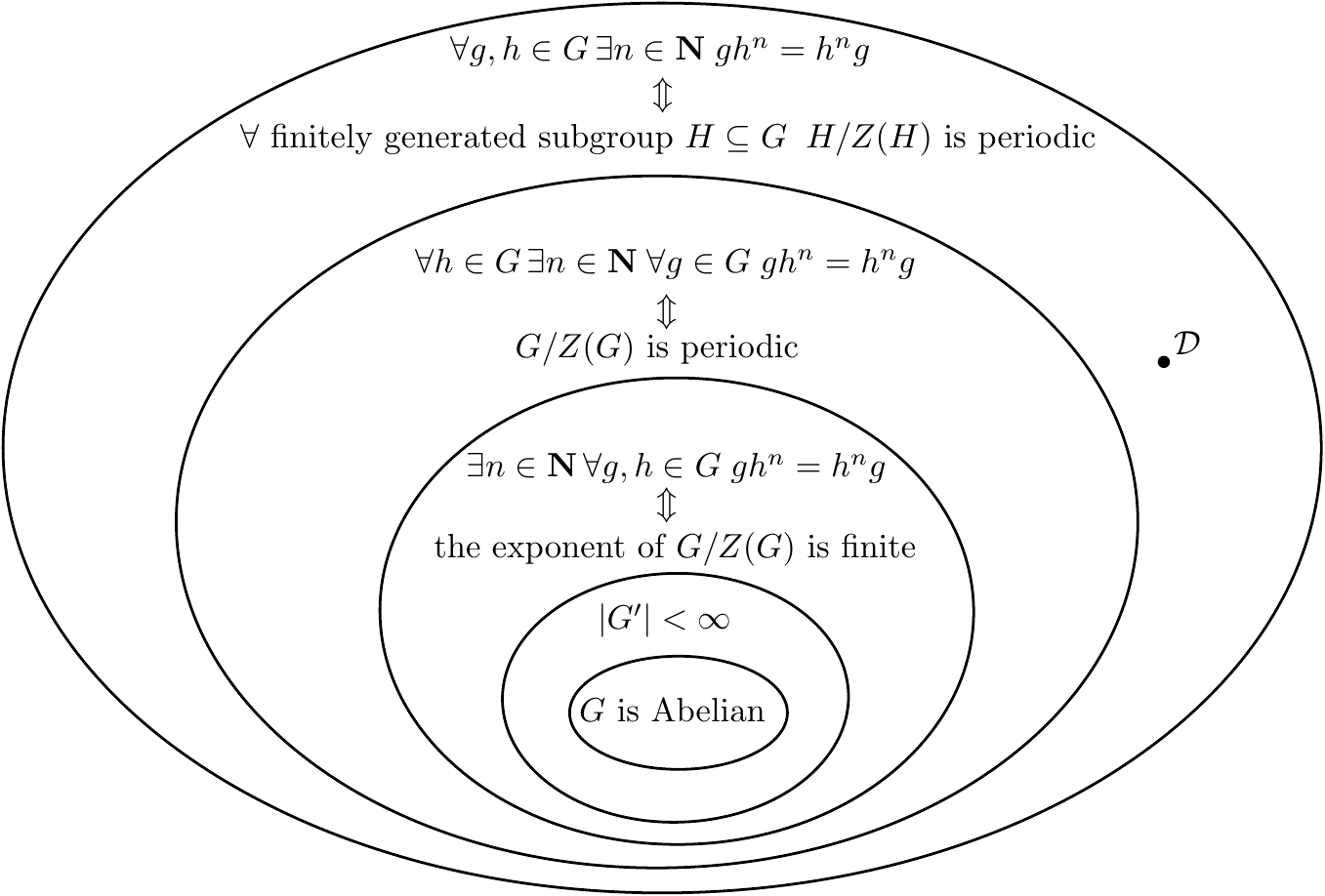}
\end{center}

\section{Proof of sufficiency in Theorems 1 and 2}

In order to prove (2) $\Rightarrow$ (1) in both theorems, like in the ungraded case, it is sufficient to find a homogeneous regular 
element in every gr-essential right ideal of given gr-(semi)prime right gr-Goldie ring $R$. Besides, we know that 
$s$ is regular iff $r_R(s)=0$ for all $s\in h(R)$.  

\vskip 5 pt
So, let $R$ be a gr-semiprime right gr-Goldie ring, $I$ a gr-essential right ideal of $R$. 

\begin{itemize}
\item {\bf Step 1.} Since the gr-Goldie dimension of $R$ is finite, one can find a gr-uniform element $a\in h(R)$ (that is, the graded right ideal $aR$ 
is gr-uniform). By \cite[Lemma 8.4.3]{N}, its right annihilator $r_R(a)$ is maximal among right annihilators of nonzero homogeneous elements of $R$.

\item {\bf Step 2.} By \cite[Lemma 8.4.2]{N}, there exists an element $b\in h(R)$ such that $ba$ is non-nilpotent. Clearly $ba$ is gr-uniform 
together with $a$. 

\item {\bf Step 3.} By induction, one can find the gr-uniform non-nilpotent elements $a_1=ba, a_2,\ldots, a_n\in I$ such that 
$a_{i+1}\in \bigcap_{j=1}^i r_R(a_j)$ for $1\leq i < n$ and $\bigcap_{j=1}^n r_R(a_j)=0$. Then, the sum $\sum_{j=1}^n a_jR$ is direct since 
$r_R(a_j)=r_R(a_j^2)$ for all $j$ (so that is why the process terminates). It follows from this that 
$r_R(a_1+\ldots +a_n)=\bigcap_{j=1}^n r_R(a_j)=0$ whence the element $a_1+\ldots +a_n\in I$ is regular. 
\end{itemize}

All these arguments are absolutely similar to the ungraded case. The main problem is that the degrees of homogeneous elements $a_1,\ldots, a_n$ 
can be different so their sum can be not homogeneous. 

\begin{proof}[\bf Finish of Proof $(2)\Rightarrow (1), (**)$ of Theorem \ref{T1}]
As noted in the introduction, this follows from the more general result of the author. We give more short proof here. 

So, suppose that $G$ is a periodic group. Let $k_i$ be an order of $\deg a_i\in G$, $i=1,\ldots, n$. Then $a_i^{k_i}\in I_e$ and 
$r_R(a_i^{k_i})=r_R(a_i)$ so $a_1^{k_1}+\ldots+a_n^{k_n}$ is a regular element from $I_e$.

Now we prove statements $(**)$. Since $rs^{-1}=(rs^{k-1})(s^k)^{-1}$, where $r\in R$, $s\in S$, $k=O(\deg s)$, we have $RS^{-1}=RS_e^{-1}$. 
Let us suppose that $p\in R_e$ is a zero divisor of $R$ and $pr=0$ $(rp=0)$ for some $r\in h(R)\setminus 0$. Then $pr^l=0$ ($r^lp=0$), 
where $l=O(\deg r)$, so $p$ is a zero divizor of $R$. Hence $S_e$ equals the set of regular elements of $R_e$. Finally, for any right ideal $I$ of 
$R_e$, we have 
$$I \text{ is essential of } R_e \Leftrightarrow IR \text{ is gr-essential of } R \Leftrightarrow IR\cap S\ne \varnothing \Leftrightarrow 
I\cap S_e=(IR\cap S)_e\ne \varnothing$$ because $R$ is $e$-faithful in view of the periodicity of $G$. Thus, $R_e$ is a semiprime Goldie ring.
\end{proof}

\begin{proof}[\bf Finish of Proof $(2)\Rightarrow (1)$ in Theorem \ref{T2}] 
Goodearl and Stafford \cite{GS} find the elements $s_2,\ldots, s_n\in h(R)$ such that $a_1^2s_2a_2^2s_3\ldots s_na_n^2\ne 0$ using 
the gr-primary of $R$. Then, they use \cite[Lemma 2]{GS} and construct the non-nilpotent elements   
$$c=s_1a_1^2s_2a_2^2s_3\ldots s_na_n^2,\quad d_i=(a_is_{i+1}a_{i+1}^2\ldots a_n^2)(s_1a_1^2\ldots s_ia_i),\;i=1,\ldots, n, \quad (s_1\in h(R))$$ 

In \cite{GS}, the group $G$ is supposed to be Abelian and so all $d_i$'s have the same degree thus their sum $d$ is homogeneous; 
the standard reasoning proves that $d$ is regular. Let us enhance this result. 

Denote $\deg c=:g$ and $\deg d_i=:h_i$, $i=1,\ldots, n$. Obviously, all $h_i$'s are conjugate. Condition (2) gives numbers 
$k_1,\ldots,k_{n-1} \in \N$ such that $h_i^{k_i}=h_{i+1}^{k_i}$ for all $i=1,\ldots, n-1$. Then $h_1^k=\ldots =h_n^k$ for $k=k_1\ldots k_{n-1}$, hence 
the homogeneous non-nilpotent elements $d_1^k,\ldots, d_n^k$ have the same degree and their sum is homogeneous. 
Since $r_R(a_i)\subseteq r_R(d_i^k)$, \cite[Lemma 3]{GS} yields 
$r_R(a_i)=r_R(d_i^k)$. Besides, the sum $\sum_{i=1}^n d_i^k R$ is direct because $d_i^kR\subseteq a_iR$. Therefore
$$r_R(d_1^k+\ldots + d_n^k)=\bigcap_{i=1}^n r_R(d_i^k)=\bigcap_{i=1}^n r_R(a_i)=0,$$ whence the element $d_1^k+\ldots +d_n^k\in h(I)$ is regular.
\end{proof}

\section{Proof of necessity in Theorems 1 and 2: the constructing of counterexamples}

\begin{proof}[\bf Proof of $(1)\Rightarrow (2)$ in Theorem 1.] 
Suppose that the group $G$ contains an element $h$ of infinite 
order. Then one can adapt the well-known example \cite[Example 9.2.2]{N2}: consider the ring $R=k[x,y]/(xy)$ ($k$ is a field) with the following 
$G$-grading: $$R_e=k, \;R_{h^n}=kx^n, \;R_{h^{-n}}=ky^n \,(n\in \N),\; R_g=0 \text{ for } g\notin \langle h\rangle.$$
We see that $Q^{gr}_{cl}(R)=R$ is gr-semiprime gr-Noetherian but not gr-Artinian ring.
\end{proof}

\begin{proof}
 [\bf Proof of $(1)\Rightarrow (2)'$ in Theorem 2.] 
Assume the converse, suppose that the group $G$ contains elements $g,h$ such that
$h^k\ne gh^lg^{-1}$ for all $k,l\in \N$. Clearly, $O(h)=\infty$.

Let us consider a polynomial ring $D=k[t]$ over field $k$ supplied the following $G$-grading: $D_{h^n}=kt^n$ for $n\geq 0$, and 
$D_\si=0$ for other $\si\in G$. Also consider a $G$-graded matrix ring
$$R=M_2(D)(e,g), \quad R_\si=\begin{pmatrix} D_\si & D_{\si g} \\ D_{g^{-1}\si} & D_{g^{-1}\si g}\end{pmatrix} \;(\si\in G).$$ 

It is clear that $R$ is a prime Noetherian (right and left) but not gr-Artinian ring.  Let us show that the set $S$ of homogeneous regular
elements of $R$ is exhausted by diagonal matrices over $k^*$. It follows from this, that $Q^{gr}_{cl}(R)=R$ and, thus, the ring $R$ 
is a counterexample to condition (1).

So, let  $A=\begin{pmatrix} a & b \\ c & d\end{pmatrix}\in S\cap R_\si$. If $b\ne 0$, then $\deg b=\si g\in \langle h\rangle$ whence
$a\in D_\si=0$, $c\in D_{g^{-1}\si}=0$ and $A$ is a zero divisor. Hence $b=0$ and similarly $c=0$. Therefore $a\ne 0\ne d$ and
$\deg d=gh^lg^{-1}=h^k$ for some $k,l\in \N_0$. In view of the choise the elements $g,h$, this is possible only if $k=l=0$, i.~e., $a,d\in k^*$.
\end{proof}

\begin{rmk} As noted in the introduction, the maximal graded right quotient ring $Q^{gr}(R)$ is gr-semisimple whenever $R$ is a gr-semiprime 
right gr-Goldie ring. Let us exemplify this for the ring $R$ from the proof $(1)\Rightarrow(2)'$. On the graded analogue of Utumi's Theorem 
\cite[proposition 2]{Kan2},
$$Q^{gr}\bigl(M_n(K)(g_1,\ldots, g_n)\bigr)=M_n(Q^{gr}(K))(g_1,\ldots, g_n)$$ for any graded ring $K$, number $n\in \N$ and set $(g_1,\ldots, g_n)\in G^n$.
Since $Q^{gr}(k[t])=k[t,t^{-1}]$ is a gr-field, we get 
$$Q^{gr}(R)=M_2(k[t,t^{-1}])(e,g)$$ is a gr-simple gr-Artinian ring.
Also note that $I=M_2(tk[t])(e,g)$ is a gr-essential ideal of $R$ without homogeneous regular elements. According to the construction of the ring 
$Q^{gr}(R)$, the matrices $t^{-1}e_{ij}\in h(Q^{gr}(R))$, where $e_{ij}$ are matrix units, $i,j\in \{1,2\}$, are the equivalence classes of 
the following homomorphisms $\phi_{ij}\colon I_R\to R_R$ or $R$-modules: 
$$\begin{aligned}
\begin{pmatrix} t^{-1} & 0 \\ 0 & 0 \end{pmatrix} \leftrightarrow \phi_{11}\colon \begin{pmatrix} ta(t) & tb(t) \\ tc(t) & td(t) \end{pmatrix}
\mapsto \begin{pmatrix} a(t) & b(t) \\ 0 & 0 \end{pmatrix}, \;\;
\begin{pmatrix} 0 & t^{-1} \\ 0 & 0 \end{pmatrix} \leftrightarrow \phi_{12}\colon \begin{pmatrix} ta(t) & tb(t) \\ tc(t) & td(t) \end{pmatrix}
\mapsto \begin{pmatrix} c(t) & d(t) \\ 0 & 0 \end{pmatrix}\!,\\
\begin{pmatrix} 0 & 0 \\ t^{-1} & 0 \end{pmatrix} \leftrightarrow \phi_{21}\colon \begin{pmatrix} ta(t) & tb(t) \\ tc(t) & td(t) \end{pmatrix}
\mapsto\begin{pmatrix} 0 & 0 \\ a(t) & b(t) \end{pmatrix}, \;\;
\begin{pmatrix} 0 & 0 \\ 0 & t^{-1} \end{pmatrix} \leftrightarrow \phi_{22}\colon \begin{pmatrix} ta(t) & tb(t) \\ tc(t) & td(t) \end{pmatrix}
\mapsto \begin{pmatrix} 0 & 0 \\ c(t) & d(t) \end{pmatrix}\!.
 \end{aligned}$$
\end{rmk}

As far as we know, the first counterexample in the gr-prime case is invented by Dmitry Bazhenov, the student of the author.

\begin{ex}[\cite{DB}] Let $k$ be a field, $G=\langle r,s\mid s^2=(rs)^2=e\rangle =\langle r\rangle_\infty \leftthreetimes \langle s\rangle_2$
an infinite dihedral group. Let us consider the ring $R=k\langle x,y,z\mid xy=yx=0,z^2=1,xz=yx,yz=zx\rangle$ with the following $G$-grading:
$$R_e=k,\quad R_{r^n}=kx^n,\quad R_{r^{-n}}=ky^n,\quad R_s=kz,\quad R_{sr^n}=kzx^n,\quad R_{sr^{-n}}=kzy^n \quad (n\in \N).$$
Since $S=k^*$ we have $Q^{gr}_{cl}(R)=R$. One can show that $R$ is gr-Noetherian but not gr-Artinian ring. 
\end{ex}

A. L. Kanunnikov found the matrix representation of the ring $R$ over a polynomial ring $k[t]$:
 $$x=\begin{pmatrix} t & 0 \\ 0 & 0 \end{pmatrix}, \quad y=\begin{pmatrix} 0 & 0 \\ 0 & t \end{pmatrix}, 
\quad z=\begin{pmatrix} 0 & 1 \\ 1 & 0 \end{pmatrix}.$$

So we can consider $R$ as a graded subring of $M_2(k[t])(e,s)$: 
$$R=\left\{ \left.\begin{pmatrix} a(t) & b(t) \\ c(t) & d(t) \end{pmatrix}\in M_2(k[t]) \,\right| a(0)=d(0),\;b(0)=c(0)\right\},$$ 

\noindent where $k[t]$ is graded in the natural way: 
$k[t]_{r^n}=kt^n$ if $n\in \N_0$ and $k[t]_g=0$ for other $g\in G$.  

Inspired by this nice example, the author generalized it to the series of rings in the proof of $(1)\Rightarrow (2)'$.

\section{Proof of $(2)'\Rightarrow (2)$ in Theorem 2}

Let $g,h$ be arbitrary elements of $G$. Show that $gh^ng^{-1}=h^n$ for some $n\in\N$. 
Consider that the order of $h$ is infinite else the statement is evident. 

On condition $(2)'$, there exist $k,l,m,n\in \N$ such that
$$hg^kh^{-1}=g^l,\qquad gh^mg^{-1}=h^n.$$
From the second equality it follows by induction that $$g^dh^{m^d}g^{-d}=h^{n^d} \eqno(\bigstar)$$
for any $d\in \N$. Indeed, to step from $d$ to $d+1$, raise $(\bigstar)$ to the $m$-th power and conjugate with respect to $g$: 
$$g^dh^{m^{d+1}}g^{-d}=h^{n^dm} \;\Longrightarrow\; g^{d+1}h^{m^{d+1}}g^{-d-1}=gh^{n^dm}g^{-1}=(gh^mg^{-1})^{n^d}=(h^n)^{n^d}=h^{n^{d+1}}.$$
  
If $k=l$, then $hg^k=g^kh$ and $h^{n^k}=g^kh^{m^k}g^{-k}=h^{m^k}$ whence $m^k=n^k$ and $m=n$. 

Consider the case $k>l$. Substituting  $h^{-1}g^lh$ for $g^k$ in $(\bigstar)$ at $d=k$, we get
$$g^kh^{m^k}g^{-k}=h^{n^k}=g^lh^{m^k}g^{-l}.$$
On the other hand, $$h^{n^k}= (h^{n^l})^{n^{k-l}}=g^lh^{m^ln^{k-l}}g^{-l}.$$ 
Comparing this with the previous equality, we have 
$$h^{m^k}=h^{m^ln^{k-l}}\; \Longrightarrow \;m^{k-l}=n^{k-l}\; \Longrightarrow \; m=n.$$

The case $k<l$ is similar.

\vskip 10 pt

\section*{Acknowledgements} 

The author thanks A. U. Olshanskii and A. A. Klyachko for helpful comments and examples of groups and D. S. Bazhenov for discussion 
that generated some ideas of this research.

\end{document}